\documentclass[12 pt]{article}  
\input amssym.def

\def\Z{{\Bbb Z}}
\def\C{{\Bbb C}}  
\def\Q{{\Bbb Q}}

\title{The Galois action on character tables}
\author{
{\sc Terry Gannon}\thanks{This is a contribution to the proceedings in honour
of John McKay. This research is supported in part by NSERC}
 \\
{\footnotesize Department of Mathematical Sciences, University of Alberta,}\\
{\footnotesize  Edmonton, Canada, T6G 2G1}\\
{\footnotesize e-mail: {\tt tgannon@math.ualberta.ca}}}
\date{}

\begin{document}
\maketitle
\rightline{To John, prophet and mentor}

\begin{abstract} 
A geometric interpretation and generalisation for the Galois action
on finite group character tables is sketched.
\end{abstract}

\section{Introduction}

Throughout this paper,
let $G$ be a finite group. Its finite-dimensional representations (over $\C$)
are completely reducible into a direct sum of irreducible representations.
Moreover, such a representation $\rho$ is uniquely determined (up to equivalence)
by its {\it character} $ch_\rho(g):={\rm Tr}(\rho(g))$. 

A character is a {\it class function}, i.e. is
constant on conjugacy classes: $ch_\rho(k^{-1}gk)=
ch_\rho(g)$. We collect the characters of irreducible representations,
evaluated at each conjugacy class, into a square matrix called
the {\it character table}. The character tables for the symmetric group $S_3$ 
and alternating group $A_5$ are given below ($\alpha=(1+\sqrt{5})/2,\alpha'=
(1-\sqrt{5})/2$). The character table is an
important group invariant, with rich structure. For example, 
each 1-dimensional representation is a symmetry permuting the rows, and
dually, any central element is a symmetry permuting the columns.   

$$\vbox{\tabskip=0pt\offinterlineskip
  \def\tablerule{\noalign{\hrule}}
  \halign to 1.5in{
    \strut#&\vrule#\tabskip=0em plus1em &    
    \hfil#&\vrule#&\hfil#&\vrule#&    
\hfil#&\vrule#&\hfil#&
\vrule#\tabskip=0pt\cr\tablerule     
&&$S_3$&&\hfil(1)\hfil&&\hfil(12)\hfil&&\hfil(123)\hfil&\cr
\tablerule&&\hfil$ch_1$\hfil&&\hfil{1}\hfil&&\hfil{1}\hfil&&\hfil{1}\hfil&\cr
\tablerule&&\hfil$ch_s$\hfil&&\hfil{1}\hfil&&\hfil${-1}$\hfil&&\hfil{1}\hfil&\cr
\tablerule&&\hfil$ch_2$\hfil&&\hfil{2}\hfil&&\hfil{0}\hfil&&\hfil${-1}$\hfil&\cr
\tablerule\noalign{\smallskip}}}
\qquad\qquad\vbox{\tabskip=0pt\offinterlineskip
  \def\tablerule{\noalign{\hrule}}
  \halign to 3in{
    \strut#&\vrule#\tabskip=0em plus1em &    
    \hfil#&\vrule#&\hfil#&\vrule#&    
\hfil#&\vrule#&\hfil#&\vrule#&\hfil#&\vrule#&\hfil#&
\vrule#\tabskip=0pt\cr\tablerule     
&&$A_5$&&\hfil(1)\hfil&&\hfil(12)(34)\hfil&&\hfil(123)\hfil&&
\hfil(12345)\hfil&&\hfil(13524)\hfil&\cr
\tablerule&&\hfil$ch_1$\hfil&&\hfil{1}\hfil&&\hfil{1}\hfil&&\hfil{1}\hfil&&\hfil{1}\hfil&&\hfil{1}\hfil&\cr
\tablerule&&\hfil$ch_4$\hfil&&\hfil{4}\hfil&&\hfil$0$\hfil&&\hfil{1}\hfil&&\hfil${-1}$\hfil&&\hfil${-1}$\hfil&\cr
\tablerule&&\hfil$ch_5$\hfil&&\hfil{5}\hfil&&\hfil{1}\hfil&&\hfil${-1}$\hfil&&\hfil{0}\hfil&&\hfil${0}$\hfil&\cr
\tablerule&&\hfil$ch_3$\hfil&&\hfil{3}\hfil&&\hfil${-1}$\hfil&&\hfil{0}\hfil&&\hfil${\alpha}$\hfil&&\hfil${\alpha'}$\hfil&\cr
\tablerule&&\hfil$ch_{3'}$\hfil&&\hfil{3}\hfil&&\hfil${-1}$\hfil&&\hfil{0}\hfil&&\hfil${\alpha'}$\hfil&&\hfil${\alpha}$\hfil&\cr
\tablerule\noalign{\smallskip}}}$$

\centerline{{\sc Table.} The character tables of $S_3$ and $A_5$}
\medskip

The Galois symmetry of the character table isn't as well-known. 
Any element of the {\it cyclotomic field} $\Q[\xi_n]$ for $\xi_n:=\exp[2\pi i/n]$
can be written (in many ways) as a polynomial $p(\xi_n)=\sum_k a_k\xi_n^k$,
$a_k\in\Q$; it is called  a cyclotomic integer if all $a_k\in\Z$. 
The entries of the character table of $G$ are cyclotomic
integers, $ch_\rho(g)\in\Z[\xi_n]$, where $n$ can 
be taken to be the order (or 
even the exponent) of $G$. This can be seen by diagonalising the matrix $\rho(g)$,
which will have order dividing $n$. 
For example in the $A_5$ table,
$\alpha=1+\xi_5^2+\xi_5^3$ and $\alpha'=1+\xi_5+
\xi_5^4$ are both cyclotomic integers.
In fact somewhat more is true: any (complex) representation of $G$ can be realised (in
many ways) by matrices with entries in $\Q[\xi_n]$ (this can be elegantly 
proved using induced representations -- see e.g. Thm. 10.3 of \cite{is}).

The Galois group ${\cal G}_n:={\rm Gal}(\Q[\xi_n]/\Q)$ is the multiplicative
group $\Z_n^\times$ consisting of all classes mod $n$ of integers coprime
to $n$. The integer $\ell\in\Z^\times_n$ corresponds to the automorphism
$\sigma_\ell\in{\cal G}_n$ taking $\xi_n$ to $\xi_n^\ell$. This implies
that $\sigma_\ell(p(\xi_n))=p(\xi_n^\ell)$ for any polynomial $p$ with
rational coefficients. 

The Galois group ${\cal G}_n\cong \Z_n^\times$ acts on the rows of the character
table of $G$, by $\sigma .ch_\rho=ch_{\sigma \rho}$, where $\sigma
\rho$ is defined by first realising $\rho$ as matrices over $\Q[\xi_n]$,
and then evaluating $(\sigma \rho)(g)=\sigma(\rho(g))$ entry-wise. That 
$\sigma\rho$ is
a representation of $G$ is clear because $\sigma$ preserves addition and
multiplication; that $\sigma\rho$ is independent of how $\rho$ is realised
over $\Q[\xi_n]$, is now clear by considering characters. For example,
$\ell\equiv \pm 2$ (mod 5) will permute the $A_5$-characters $ch_3$ and
$ch_3'$. (Of course, $\sigma\rho$ can also be expressed in a basis-independent way
e.g. by having $\sigma\in{\rm Gal}(\C/\Q)$ act on $\rho\in{\rm Hom}(G,{\rm End}
(\C^m))$ by $\sigma\circ\rho\circ\sigma^{-1}$.)

The Galois group $\Z^\times_n$ acts on columns of character tables by
$\sigma_\ell.g=g^\ell$ ($\ell$ is only defined mod $n$, but $g^n=1$ so
this isn't a problem). That this takes conjugacy classes to conjugacy classes
is clear: $(k^{-1}gk)^\ell=k^{-1}g^\ell k$.

These two actions are compatible, in the sense that $\sigma_\ell$ applied
to the number $ch_\rho(g)\in \Q[\xi_n]$ is given by both the row- and 
column-actions:
\begin{equation}\label{compat}
\sigma_\ell(ch_\rho(g))=ch_{\sigma_\ell \rho}(g)=ch_\rho(g^\ell)\ .
\end{equation}
The first equality is clear by the definition of $\sigma\rho$: Tr$(\sigma
(\rho(g)))=\sigma({\rm Tr}(\rho(g)))$. The second equality follows by
diagonalising $\rho(g)$ (which can be done over $\Q[\xi_n]$): its eigenvalues
are $n$th roots of unity so will be raised to the $\ell$th power.

For example, the $S_3$-character values are all integers, and so are 
invariant under any Galois automorphism. This is equivalent to saying that
$g^\ell$ is conjugate to $g$, for any $g\in S_3$ and any $\ell\in\Z_6^\times$.
Since $g^\ell$ will likewise be conjugate to $g$ in any symmetric group
$S_m$, all character values for any $S_m$ will necessarily be integers.

There are many theorems in group theory whose best (or only) known proof uses
character table technology, including this Galois action. A famous example
is Burnside's $p^aq^b$-theorem, which says that as long as three distinct primes 
don't divide the order of a group, that group will necessarily be solvable. 
This is the best possible result, in the sense that the alternating group
$A_5$ is not solvable (much to the chagrin of quintic solvers) and only
three primes divide $|A_5|=60$.

We can formulate this Galois action more abstractly as follows. The characters
are class functions, so lie in (in fact span) the space
\begin{equation}
{\rm  Map}_G(G,\Q[\xi_n])\label{dum}
\end{equation}
of $G$-equivariant maps, where $G$ acts on itself by conjugation and fixes the 
scalars. The Galois group
${\cal G}_n$ acts on this space in two ways: it manifestly acts on the 
scalars, and it acts on $G$ as above. The compatibility (\ref{compat})
of these two actions
is clearly a basis-dependent statement; that such a basis (e.g. 
the irreducible characters) exists, is equivalent to the statement that
the representation of ${\cal G}_n$ on the group algebra $\Q \langle g\rangle$,
for any $g\in G$, is a subrepresentation of that on the field $\Q[\xi_n]$.

In any case, it can be asked why should the cyclotomic Galois group ${\cal G}_n$
act on finite group character tables. More precisely, these Galois groups
act by definition on cyclotomic numbers; any other action that can't be
reduced quite directly to this one can be termed `mysterious'.
In this sense, perhaps the action of ${\cal G}_n$ on $G$ can be labelled
mysterious. Can we relate it to more familiar Galois actions, perhaps as
 a special case of some general class? 
I think this is a question John McKay could ask. 
We'll suggest an answer in the following pages. 

However, especially in this era of George Bush and Tony Blair, we've been 
taught that we shouldn't look too closely at motivations. 
Indeed, for an additional ulterior motive for this material, see Chapter 6 
of \cite{gan}.  
Certainly, the `explanation' we propose in Section 3 for this action is much 
more complicated than the standard proof sketched above. At the end of the day,
we should judge a question by how intriguing the resulting picture is, which 
it helped shape. So withhold your judgement until you conclude Section 3!

\section{A warm-up example}

To help clarify what is wanted here, let's consider a toy model. It is a
vaguely similar action, though is much simpler to `explain'. What we are after
is essentially its profinite version. This section can be skipped if
the material is too unfamiliar; on the other hand, it may also be of
independent interest.

In Moonshine and orbifold string theory, the following action of SL$_2(\Z)$ arises.
Let $G$ be any finite group. Consider all pairs $(g,h)\in G\times G$ of
commuting pairs. Identify simultaneous conjugates: 
$(g,h)\sim(k^{-1}gk,k^{-1}hk)$. For example, $S_3$ has 8 such equivalence
 classes, while $A_5$ has 22. On these equivalence classes of pairs, we 
get a (right) SL$_2(\Z)$ action by
\begin{equation}\label{modact}
(g,h)\left(\matrix{a&b\cr c&d}\right)=(g^ah^c,g^bh^d)\ .
\end{equation}

Where does this come from? The answer is geometry. If we forget the commuting
pairs and simultaneous conjugation, the SL$_2(\Z)$-action becomes a braid
group ${\cal B}_3$-action:
\begin{equation}\label{bract}
(g,h).\sigma_1=(g,gh)\ ,\qquad(g,h).\sigma_2=(gh^{-1},h)\ ,
\end{equation}
where ${\cal B}_3=\langle \sigma_1,\sigma_2\,|\,\sigma_1\sigma_2\sigma_1=
\sigma_2\sigma_1\sigma_2\rangle$.

Now, one of the basic facts about the braid groups \cite{bir} is that they act on the
free groups. In particular, ${\cal B}_3$ acts on ${\cal F}_3=\langle x,y,z\rangle$ by
$$\sigma_1.x=xyx^{-1}\ ,\ \sigma_1.y=x\ ,\ \sigma_1.z=z\ ,\ \sigma_2.x=x\ ,\
\sigma_2.y=yzy^{-1}\ ,\ \sigma_2.z=y\ .$$ 
The easiest way to see this is through the realisation of the braid group
as the fundamental group of configuration space. Indeed, the pure braid
group ${\cal P}_n$ is
$\pi_1(C_n)$ where $C_n=\{(z_1,z_2,\ldots,z_n)\in\C^n\,|\,z_i\ne z_j\})$
and the free group ${\cal F}_n$ is $\pi_1(\C-\{1,2,\ldots,n\})$, so given the obvious
 bundle $C_{n+1}\rightarrow C_n$ (dropping a coordinate), there will be an action of the fundamental
group ${\cal P}_n$ of the base on the fundamental group ${\cal F}_n$ of
the fibre. The ${\cal B}_n$-action is obtained through symmetrisation of $C_n$.

In any case, ${\cal B}_3$ acts on ${\cal F}_3$ (by group automorphisms), so it 
also acts on the {\it set} Hom$({\cal F}_3,
G)\cong G^3$ (this latter action won't respect the $G\times G\times G$ group structure).
Collapsing a triple $(g_1,g_2,g_3)\in G^3$ to the pair $(g_1g_2^{-1},g_2g_3^{-1})$
recovers (\ref{bract}). This ${\cal B}_3$-action becomes a well-defined SL$_2(\Z)$-action
if we either restrict to a commuting pair, or we identify simultaneous
conjugates. 

Similarly, we get a ${\cal B}_n$-action on $G^n\cong {\rm Hom}({\cal F}_n,G)$
and hence on $G^{n-1}$. The role of SL$_2(\Z)$ will now be served by
${\cal B}_n/\langle z\rangle$ ($n$ even) or ${\cal B}_n/\langle z^2\rangle$
($n$ odd), where $z$ is the generator of the centre of ${\cal B}_n$.
This `explains' the SL$_2(\Z)$-action (\ref{modact}) in the sense that it embeds it
in a larger context, showing that it belongs to an infinite family of related
actions, hence generalising it.
We want to do the same to the Galois action on character tables -- to embed
it in an infinite family of other Galois actions. In this 
 geometric interpretation of Galois on character tables,
${\cal B}_n$ is replaced by Gal$(\overline{\Q}/\Q)$ and ${\cal F}_n$ by
its profinite completion.

Of course there is never a unique `explanation', nor need there even be a best one.
For a different explanation, (\ref{modact}) can also be obtained from the following construction
(see e.g. \cite{tur,gan} for the background). Given a Hopf algebra, the Drinfel'd double produces
a quasi-triangular Hopf algebra (i.e. a Hopf algebra co-commutative up to
an isomorphism). Given such a quasi-triangular Hopf algebra, Reshetikhin--Turaev
obtain from this a modular category. These are interesting because modular
categories supply knot invariants on arbitrary 3-manifolds. In a modular
category (or what is essentially the same thing, a 3-dimensional topological
field theory or 2-dimensional rational conformal field theory), we get
representations of any surface mapping class group. Now apply this construction
to the group algebra $\C G$ (a Hopf algebra). Fix any compact genus $g$
surface $\Sigma_g$, and consider the set of all group homomorphisms $\pi_1(\Sigma_g)\rightarrow G$,
where we identify two maps if they are conjugates by $G$.
The mapping class group $\Gamma_{g,0}$ acts naturally on $\pi_1(\Sigma_g)$
and hence on Hom$(\pi_1(\Sigma_g),G)$, and in genus $g=1$ this
action reduces to that of (\ref{modact}) (recall that $\pi_1$ of a torus is $\Z^2$
and $\Gamma_{1,0}$ is SL$_2(\Z)$).

The second explanation certainly should be regarded as an interesting 
generalisation of (\ref{modact}). Perhaps one advantage though of the first
explanation is that it is both more intrinsically group-theoretic and more
elementary, not relying on relatively
heavy machinery. In any case, the interpretation
we'll suggest shortly for the Galois action on character tables is
related to this first explanation.

\section{The geometric meaning of Galois}

Seeking a geometric interpretation of a Galois action leads naturally to
geometric Galois. Classical Galois theory concerns field extensions $L$
over a base field $K$; geometric Galois concerns unramified coverings $f:Y\rightarrow X$
of a space $X$ (see e.g. \cite{kug} for a gentle account). The role of the  algebraic
closure $\overline{K}$ is played by the universal covering space $\widetilde{X}$,
and the Galois group Gal$(\overline{K}/L)$ by the fundamental group $\pi_1(Y)$.
The action of Galois automorphisms on $\alpha\in L$ corresponds to that of
deck transformations on the fibre $f^{-1}(p)$ above $p\in X$ (a {\it deck} or {\it covering 
transformation} is a homeomorphism $\gamma:Y\rightarrow Y$ satisfying
$f\circ\gamma=\gamma$; we identify the deck transformation with the permutation
it induces on the fibre $f^{-1}(p)$). The deck transformations 
 form a group analogous to Gal($L/K)$). 
A {\it Galois covering} is one whose deck transformation group is transitive on
$f^{-1}(p)$; if the degree $n$ of the covering is finite, this is the same
as saying the order of that group is also $n$.

In complete analogy with classical Galois theory, geometric Galois establishes
a bijection between subgroups $H$ of $\pi_1(X)$ and coverings $Y=\widetilde{X}/H$,
where $H\cong \pi_1(Y)$. $Y$ is a Galois covering iff $H$ is normal in $\pi_1(X)$,
in which case the deck transformation group is isomorphic to
$\pi_1(X)/H$. 

For example, consider the covering (spiral staircase) $Y_n: x_1^nx_2=1$ of
the punctured complex plane $X:x_1x_2=1$. Here, $\pi_1(X)\cong \Z$ and $\pi_1(Y_n)$
can be naturally identified with the subgroup $n\Z$. The deck transformations
are $(x_1,x_2)\mapsto (\xi_n^kx_1,x_2)$, and form the group $\Z_n\cong\Z/n\Z$.

Geometric Galois would have been understood shortly after Poincar\'e defined
the fundamental group -- it can be understood in terms of function fields. 
But it took Grothendieck much later to explain how to
simultaneously combine geometric Galois with the Galois theory of number fields. 
Let $X$ be an
algebraic variety which can be defined over $\Q$, i.e. it's given by a set of 
polynomials $p_i(x_1,\ldots,x_n)$
with rational coefficients. Fix (if it exists) a rational basepoint $p\in X(\Q)$.
Choose any finite index normal subgroup $N$ of the topological fundamental group
$\pi_1(X)$. Then there will be
a covering $f:X_N\rightarrow X $ with $\pi_1(X_N)=N$. The Generalised Riemann Existence
Theorem ensures that $X_N$ will be an algebraic variety over the algebraic
closure $\overline{\Q}$ (i.e. the coefficients of the defining polynomials
will be algebraic numbers). In particular, the fibre $f^{-1}(p)$ will all
have coordinates in $\overline{\Q}$.

The absolute Galois group $\Gamma_{\Q}={\rm Gal}(\overline{\Q}/\Q)$ acts
on the points in $f^{-1}(p)$ component-wise. It acts on varieties over $\overline{\Q}$, equivalently the finite index subgroups $N$, coefficientwise on the
defining polynomials. Grothendieck wrote that an automorphism 
$\sigma\in\Gamma_{\Q}$ acts on the deck transformations
by
\begin{equation}\label{groth}
\gamma_N\mapsto \sigma\circ \gamma_{\sigma^{-1}N}\circ\sigma^{-1}\ ,
\end{equation}
where as always we identify a deck transformation with the associated permutation
of $f^{-1}(p)$. 

This isn't quite well-defined yet because $\sigma$ may mix up
the deck transformations for different coverings -- 
 when $\sigma N\ne N$ we must explain how 
$\gamma_{\sigma^{-1}N}$ is related to $\gamma_N$. This isn't hard, because
for each $\sigma$ the deck transformation groups $\pi_1(X)/(\sigma N)$ are 
all naturally isomorphic
(the isomorphism being given by the appropriate Galois automorphism). A more 
elegant solution to this well-definedness issue will be given shortly.
 But an example should help clarify
the basic idea of (\ref{groth}). Return to the covering $Y_n$ of the punctured plane $X$.
Take $p=(1,1)$, so $f^{-1}(p)=\{(\xi_n^j,1)\}$. The deck transformation
$\gamma_k$, for $k\in\Z_n$, sends $(\xi_n^j,1)$ to $(\xi_n^{j+k},1)$. Since $Y_n$
is also defined over $\Q$, any $\sigma\in\Gamma_{\Q}$  fixes it (hence
also fixes $N=n\Z$). Note also that the $\Gamma_\Q$-action here factors through
to that of the cyclotomic group ${\cal G}_n\cong\Z_n^\times$: then (\ref{groth})
says that $\ell\in\Z^\times_n$
acts on the deck transformation $\gamma_k$ by
$$(\xi_n^j,1)\mapsto(\xi_n^{\ell(\ell^{-1}j+k)},1)=(\xi_n^{j+\ell k},1)\ ,$$
i.e. $\Z^\times_n$ acts on $\Z_n$ by multiplication.

There is a cleaner way to formulate this, which in addition makes 
(\ref{groth}) well-defined. All of these normal subgroups $N$ can be handled
simultaneously, by using the inverse-limit (algebra's way to integrate):
the inverse-limit lim$_\leftarrow$ of the deck transformation groups
 $\pi_1(X)/N$ is called the {\it profinite completion} $\widehat{\pi_1(X)}$
of $\pi_1(X)$ and the {\it algebraic fundamental group} of $X$.
Grothendieck's (\ref{groth}) defines `component-wise' the $\Gamma_\Q$-action on $\widehat{\pi_1(X)}$.

For example, return to the coverings $Y_n$ of the punctured plane $X$. The
algebraic fundamental group here is lim$_\leftarrow \Z_n=\widehat{\Z}$,
the direct product over all primes $p$ of the $p$-adic integers $\widehat{\Z}_p$.
It can be realised as the collection of all sequences $\hat{k}=(k_1,k_2,\ldots)\in
\prod_{n=1}^\infty\Z_n$ with the property that $k_n\equiv k_{m}$ (mod $n$)
whenever $n$ divides $m$. (Any other profinite completion is constructed
similarly.) Then $\sigma\in\Gamma_\Q$ acts on $\widehat{\Z}$ through the
{\it cyclotomic character} $\chi_{{cyclo}}$ defining the projection of $\Gamma_\Q$
onto the abelianisation $\Gamma_\Q/[\Gamma_\Q,\Gamma_\Q]$. By the Kronecker--Weber
theorem, the abelianisation is simply Gal$(\cup_n\Q[\xi_n]/\Q)={\rm lim}_\leftarrow
{\cal G}_n=\widehat{\Z}^\times$, the invertible elements in $\widehat{\Z}$,
i.e. all sequences $\hat{\ell}=(\ell_1,\ell_2,\ldots)\in\widehat{\Z}$ such that each
$\ell_n$ is coprime to $n$. The action of $\sigma\in\Gamma_\Q$ on $\widehat{\Z}$
is given by multiplication by $\chi_{cyclo}(\sigma)$.

Section 2 tells us how to proceed. We get an action of $\widehat{\Z}^\times$,
or if you prefer $\Gamma_\Q$, on the (continuous) homomorphisms Hom($\widehat{\Z},G)\cong G$, namely $\hat{\ell}$ sends $g\in G$ to $g^{\ell_n}$ where $n$ is
the order of $g$. This is precisely the Galois action on columns of the character table!
The more straightforward action on rows can be recovered from this
by a duality argument.

In other words, we should replace the definition (\ref{dum}) of the space of
class functions, with
\begin{equation}
{\rm Map}_G({\rm Hom}(\widehat{\Z},G),\overline{\Q})\ .\label{better}
\end{equation}
This expression is superior to (\ref{dum}) in that 
the two Galois actions on the class functions are both manifest in (\ref{better}):
$\Gamma_{\Q}$ acts on $\widehat{\Z}$ much as it acts on $\overline{\Q}$.

But why should we limit ourselves to $X$ being the punctured plane? Take $X$
for instance to be
$\C$ with $n$ punctures. This can always be taken to be an algebraic variety 
over $\Q$ --
e.g. the polynomial $z_1z_2^2(z_2-1)^2=z_2^2-z_2+1$ realises the plane with
two punctures. For $n$ punctures, Grothendieck's action (\ref{groth}) 
becomes a homomorphism $Gr_G^n:\Gamma_\Q\rightarrow{\rm Aut}(G)$, for any
finite group $G$ with $n$ generators. It can be reclothed into an action
of $\Gamma_{\Q}$ on $\widehat{{\cal F}_n}$ (again by group automorphisms), and
hence an action of $\Gamma_\Q$ on the {\it set} $G^n$, for any finite group 
$G$.  A minor
miracle is that, by Belyi's Theorem (see e.g. \cite{sch}), these are all 
faithful actions for $n\ge 2$ (as $G$ varies)! This faithfulness is 
very surprising, considering how large the kernel is for $n=1$. 

The homomorphism $Gr^n_G$ is a generalisation to $n>1$ of the fact that 
$g\mapsto g^\ell$ is an automorphism of cyclic $G$ (but not of general $G$).
This $\Gamma_{\Q}$-action {\it stabilises}, in the sense that $Gr^n_G=Gr^m_G$
whenever $G$ has at most $m$ and $n$ generators, so we can unambiguously 
write $Gr_G$.

The $\Gamma_{\Q}$-action on $G^n$  can be described more explicitly as follows:
\begin{equation}
\sigma.(g_1,\ldots,g_n)=(Gr_H(g_1),\ldots,Gr_H(g_n))\end{equation}
 where $H=\langle
g_1,\ldots,g_n\rangle$. In particular, $\sigma$ commutes both with 
simultaneous conjugation and with permutations of the components: 
$\sigma.(k^{-1}g_1k,
\ldots,k^{-1}g_nk)=k^{-1}(\sigma.(g_1,\ldots,g_n))k$ for any $k\in G$,
and $\pi(\sigma.(g_1,
\ldots,g_n))=\sigma.(g_{\pi 1},\ldots,g_{\pi n})$ for any permutation $\pi\in
{\cal S}_n$. Therefore (\ref{better}) generalises to
\begin{equation}
{\rm Map}_G({\rm Hom}(\widehat{{\cal F}_n},G),\overline{\Q})/{\cal S}_n\ ,
\label{best}
\end{equation}
where the symmetric group ${\cal S}_n$ acts by permuting the $n$-tuples.
Again we have two manifest (and for $n\ge 2$, faithful!)  $\Gamma_\Q$-actions. 

The most obvious question is, are there any group-theoretic consequences of 
this
action, when $n\ge 2$? Is there perhaps some generalisation of character tables
which realise this action, i.e. some generalisation of character which can
supply a natural basis for the space (\ref{best})? The obvious answer is very
classical: the $n$-characters (essentially
the coefficients of Frobenius' group determinant \cite{jo}) 
span the space (\ref{best}) (with a little help). However, these are always 
$\Q[\xi_{|G|}]$-valued, so
compatibility (\ref{compat}) will fail for them for $n>1$. In this sense
the $n$-characters are too close to the usual characters; there should be (from
this Galois perspective) a better basis of (\ref{best}).

There are plenty of other more recent multi-variable generalisations
of character tables. For instance, note that the subspace of $n$-tuples
$(g_1,\ldots,g_n)$ of mutually commuting elements will be mapped to itself by
$\Gamma_\Q$; this suggests considering the spaces of Section 2, or the 
generalisations of characters associated with elliptic 
cohomology (see e.g. \cite{kuh}). It may also be interesting to study
this Galois action specifically for the reductive groups over finite fields,
where we have a more geometric alternative to the character table, 
namely character sheaves \cite{Lu} (I thank Clifton Cunningham for this
interesting suggestion).

The case of the twice-punctured plane is especially accessible, because of
its connection to dessins d'enfants \cite{sch} and also to modular curves \cite{scho}: its coverings are
precisely the algebraic curves defined over $\overline{\Q}$, and they
are isomorphic to $({\Bbb H}\cup\Q\cup\{i\infty\})/\Gamma$ for finite index
subgroups of SL$_2(\Z)$.  For example, explicit genus-0 calculations can be 
found in \cite{cg}.

For example, Galois coverings with deck transformation
group isomorphic to the dihedral group $D_n$ (of order $2n$) are given by
$f(x)=(2-x^n-x^{-n})/4$. The fibre $f^{-1}(1/2)$ are the $2n$ points $x\in
\{\xi_{4n}^j\,|\,j=\pm 1,\pm 3,\ldots,\pm(2n-1)\}$. The deck transformations
are parametrised by pairs $(\epsilon,k)$ where $k\in\Z_n$ and $\epsilon=\pm 1$;
it acts on $f^{-1}(1/2)$ by $(\epsilon,k).\xi_{4n}^j=\xi_{4n}^{\epsilon j+4k}$.
These coverings are also defined over $\Q$, and the Grothendieck action (\ref{groth})
again factors through the cyclotomic character: $Gr_{D_n}(\sigma)$ sends 
$(\epsilon,k)$ to $(\epsilon,\chi_{cyclo}(\sigma).k)$ where $\hat{\ell}.k=\ell_nk$.
Note that this $\widehat{\Z}^\times$-action differs from the $\widehat{\Z}$-action
$g\mapsto g^{\hat{\ell}}$ coming from the once-punctured plane: for the
latter, $(-1,k)$ would be fixed by any $\hat{\ell}$, since it has order 2. 
Nevertheless the
dihedral groups are close enough to the cyclic ones that the cyclotomic
character arises again. However as mentioned earlier, for any $\sigma\in\Gamma_{\Q}$, 
there will be a finite group $G$ (with two generators) such that $\sigma$
acts nontrivially on $G$.

This short paper has tried to identify a McKay-esque question: where 
does the $g\mapsto g^\ell$ symmetry of character tables really come from? 
We sketched a geometric source for it,
relating it to more familiar Galois actions, thus fitting the $g\mapsto
g^\ell$ action into the first window of an infinite tower. We propose
that these form a natural generalisation of $g\mapsto g^\ell$.
The big question
now is, do these other Galois actions on finite groups play any role in group 
theory? This question, and several others begged by this interpretation,
 will be addressed in future work.


\begin{thebibliography}{99}

\bibitem{bir} J. S. Birman, {\it Braids, Links, and Mapping Class Groups} (Princeton
University Press, 1974).

\bibitem{cg} J.-M. Couveignes and L. Granboulan, ``Dessins from a geometric point
of view'', in: {\it The Grothendieck Theory of Dessins d'enfant}, (Cambridge University
Press, 1994) 79--113.

\bibitem{gan} T. Gannon, {\it Moonshine Beyond the Monster} (Cambridge University
Press, 2006).


\bibitem{is} I. M. Isaacs, {\it Character Theory of Finite Groups} (Dover, Nineola N.Y.
1994).

\bibitem{jo} K. W. Johnson, ``The Dedekind--Frobenius group determinant: New life
in an old problem'', {\it Groups St. Andrews (Bath, 1997),} Vol. II (Cambridge
University Press, 1999) 417--428.

\bibitem{kug} M. Kuga, {\it Galois' Dream: Group theory and differential equations}
 (Birkh\"auser, Boston 1993).

\bibitem{kuh} N. J. Kuhn, ``Character rings in algebraic topology'', in:
{\it Advances in Homotopy Theory} (Cambridge University Press, 1989) 111--126.

\bibitem{Lu} G. Lusztig, ``Character sheaves and generalizations'', in:
{\it The Unity of Mathematics} (Birkh\"auser, Boston 2006) 443--455.

\bibitem{sch} L. Schneps, ``Dessins d'enfants on the Riemann sphere'', in: 
{\it The Grothendieck Theory of Dessins d'enfant} (Cambridge University Press, 1994) 47--77.

\bibitem{scho} A. J. Scholl, ``Modular forms on noncongruence subgroups'', in:
{\it S\'eminaire de Th\'eorie des Nombres, Paris 1985-86} (Birkh\"auser, Boston 1987) 199--206.

\bibitem{tur} V. G. Turaev, {\it Quantum Invariants of Knots and 3-manifolds}
(de Gruyter, Berlin 1994).

\end{thebibliography}
\end{document}